%
\magnification 1200
\input amssym.def

\def\ifundefined#1{\expandafter\ifx\csname#1\endcsname\relax}

\newif\ifdevelop \developfalse

\newtoks\sectionnumber
\newcount\equationnumber

\nopagenumbers
\newtoks\runninghead
\headline={\ifdevelop{\tenrm\number\day--\the\month--\number\year
 \hss\the\runninghead\hss
 p.\if-\the\sectionnumber\else\the\sectionnumber.\fi\folio}\else
 \ifnum\pageno=1 \hfill \else\hss{\tenrm--\folio--}\hss \fi\fi}
\def\assignnumber#1#2{%
 \ifundefined{#1}\relax\else\message{#1 already defined}\fi
 \expandafter\xdef\csname#1\endcsname
 {\if-\the\sectionnumber\else\the\sectionnumber.\fi\the#2}}%

\let\sPP=\smallbreak
\let\mPP=\medbreak
\let\bPP=\bigbreak

\def\sLP{\smallbreak\noindent}

\def\bLP{\bigbreak\noindent}

\def\Sec#1 #2 {\ifdevelop\vfill\eject\pageno1\else
	\vskip0pt plus.1\vsize\penalty-250\vskip0pt plus-.1\vsize
	\bigbreak\bigskip\fi
	\sectionnumber{#1} \equationnumber0
 \noindent{\bf\if-#1\else #1. \fi#2}\par
	\nobreak\smallskip\noindent}

\def\wrt{with respect to }

\def\al{\alpha}
\def\be{\beta}
\def\ga{\gamma}
\def\de{\delta}

\def\la{\lambda}
\def\si{\sigma}

\def\iy{\infty}

\def\const{{\rm const.}\,}

\def\half{{\scriptstyle{1\over2}}}
\def\thalf{{\textstyle{1\over2}}}

\def\eq#1{\relax
	\global\advance\equationnumber by 1
	\assignnumber{EN#1}\equationnumber
	\ifdevelop{\rm(#1)}%
	\else{\rm (\csname EN#1\endcsname)}\fi%
	}

\def\eqtag#1{\relax\ifdevelop(#1)\else\eqtagelse{#1}\fi}

\def\eqtagelse#1{\ifundefined{EN#1}\message{EN#1 undefined}{\sl (#1)}%
	\else(\csname EN#1\endcsname)\fi%
	}

\font\titlefont=cmr10 scaled \magstep1
\font\addressfont=cmr10 at 10truept
\font\ttaddressfont=cmtt10 at 10truept

\centerline
{\titlefont
Uniform multi-parameter limit transitions in the Askey tableau}
\bPP
\centerline{Tom H. Koornwinder}
\bLP
\font\smallfont=cmr10 at 10 truept
{\smallfont Extended abstract for the Proceedings of the Conference
``Modern developments in complex analysis and related topics''
(on the occasion of the 70th birthday of prof.\ dr.\ J. Korevaar),
University of Amsterdam, January 27--29, 1993.}

\Sec1 {Elementary limit formulas}
We consider the classical orthogonal polynomials as monic polynomials
$p_n(x)=x^n+$ terms of degree less than $n$. We have
\item{$\bullet$}
Jacobi polynomials $p_n^{(\al,\be)}(x)$ \wrt weight function
$(1-x)^\al\,(1+x)^\be$ on $(-1,1)$;
\item{$\bullet$}
Laguerre polynomials $\ell_n^\al(x)$ \wrt weight function $e^{-x}\,x^\al$
on $(0,\iy)$;
\item{$\bullet$}
Hermite polynomials $h_n(x)$ \wrt weight function $e^{-x^2}$ on $(-\iy,\iy)$.

\sLP
Note that, for $\al\to\iy$, the rescaled Jacobi weight function
$(1-x^2/\al)^\al$ on $(-\al^\half,\al^\half)$ tends to the Hermite
weight function $e^{-x^2}$ on $(-\iy,\iy)$.
Accordingly we have the limit formula
$$
\lim_{\al\to\iy}\al^{n/2}p_n^{(\al,\al)}(x/\al^\half)=h_n(x).
$$
Also, for $\be\to\iy$, the rescaled Jacobi weight function
$x^\al\,(1-x/\be)^\be$ on $(0,\be)$ tends to the Laguerre weight function
$x^\al e^{-x}$ on $(0,\iy)$.
Accordingly we have the limit formula
$$
\lim_{\be\to\iy}(-\be/2)^n\,p_n^{(\al,\be)}(1-2x/\be)=\ell_n^\al(x).
$$
This can be graphically indicated in the $(\al,\be)$-parameter plane
extended with the lines $\{(\al,\be)\mid\al=\iy,\;-1<\be\le\iy\}$ and
$\{(\al,\be)\mid-1<\al\le\iy,\; \be=\iy\}$.
When we start with a point $(\al,\al)$ then we can draw a diagonal arrow
to the (Hermite) point $(\iy,\iy)$ and a vertical arrow to the (Laguerre)
point $(\al,\iy)$.

The celebrated Favard theorem states that $\{p_n\}_{n=0,1,2,\ldots}$
is a system of monic orthogonal polynomials \wrt a positive orthogonality
measure if and only if a recurrence relation
$$
\eqalignno{
x\,p_n(x)&=p_{n+1}(x)+B_n\,p_n(x)+C_n\,p_{n-1}(x),
\quad n=1,2,\ldots,&\eq{3}
\cr
x\,p_0(x)&=p_1(x)+B_0\,p_0(x),
\cr
p_0(x)&=1,
\cr}
$$
is valid with $C_n>0$ and $B_n$ real.
Below, when we will give this recurrence relation with explicit
coefficients then we will silently assume that the case $n=0$ has the
same analytic form as the case $n>0$, but with the term
$C_n\,p_{n-1}(x)$ omitted.

If the coefficients $B_n$ and $C_n$ are given then $p_n$ is completely
determined by this recurrence relation. In particular, if $B_n$ and $C_n$
would continuously depend on some parameter $\la$ then $p_n(x)$ will also
continuously depend on $\la$.
For example, Hermite polynomials satisfy the recurrence relation
$$
x\,h_n(x)=h_{n+1}(x)+\thalf n\,h_{n-1}(x).
\eqno\eq{1}
$$
Now consider rescaled Laguerre polynomials
$$
p_n(x)=p_n(x;\al,\rho,\si):=\rho^n\,\ell_n^\al(\rho^{-1}x-\si).
$$
{}From the well-known recurrence relation for Laguerre polynomials
we find for these rescaled polynomials:
$$
x\,p_n(x)=p_{n+1}(x)-\rho\,(2n+\al+1+\si)\,p_n(x)+
\rho^2\,n\,(n+\al)\,p_{n-1}(x).
\eqno\eq{2}
$$
We would like to make the rescaling in such a way that, as $\al\to\iy$,
$p_n(x)$ will tend to $h_n(x)$.
It is easy to see how to do this when we compare \eqtag{1} and \eqtag{2}.
Put $\rho:=(2\al)^{-\half}$, $\si:=-\al$. Then \eqtag{2} becomes
$$
p_{n+1}(x)-(2\al)^{-\half}\,(2n+1)\,p_n(x)+{n\,(n+\al)\over 2\al}\,p_{n-1}(x).
$$
The recurrence coefficients now tend to 0 resp.\ $\thalf n$ as $\al\to\iy$.
Hence $p_n(x)\to h_n(x)$ as $\al\to\iy$, i.e.,
$$
\lim_{\al\to\iy} (2\al)^{-\half n}\,\ell_n^\al((2\al)^\half x+\al)=h_n(x).
$$

Thus, in the extended $(\al,\be)$-parameter plane we can also start
at a Laguerre point $(\al,\iy)$ and draw a horizontal arrow to
the Hermite point $(\iy,\iy)$.

\Sec2 {Uniform limit of Jacobi polynomials}
It is now natural to conjecture that we might also make these limit transitions
in the parameter plane in a more uniform way, i.e., to make such a rescaling
of the Jacobi polynomials that they depend continuously on $(\al,\be)$ in the
extended parameter plane and reduce to (possibly rescaled) Laguerre and
Hermite polynomials on the boundary lines and boundary vertex at infinity,
respectively.
For this purpose we consider Jacobi polynomials with arbitrary rescaling:
$$
p_n(x):=\rho^n\,p_n^{(\al,\be)}(\rho^{-1}x-\si).
\eqno\eq{6}
$$
These polynomials satisfy recurrence relations \eqtag{3} with
$$
\eqalign{
C_n:=&
\rho^2\,{4n\,(n+\al)\,(n+\be)\,(n+\al+\be)\over
(2n+\al+\be-1)\,(2n+\al+\be)^2\,(2n+\al+\be+1)}
={\rho^2\,\al\,\be\,(\al+\be)\over(\al+\be)^4}
\cr
&\quad\times{4n\,(1+n/\al)\,(1+n/\be)\,(1+n/(\al+\be))\over
(1+(2n-1)/(\al+\be))\,(1+2n/(\al+\be))^2\,(1+(2n+1)/(\al+\be))}\
\cr}
\eqno\eq{4}
$$
and
$$
\eqalignno{
B_n:=&\rho\,\Bigl({\be^2-\al^2\over(2n+\al+\be)\,(2n+\al+\be+2)}+\si\Bigr)
\cr
=&\rho\,\Bigl({\be-\al\over\be+\al}\,
{1\over1+2n/(\al+\be)}\,{1\over1+(2n+2)/(\al+\be)}+\si\Bigr).&\eq{5}
\cr}
$$
{}From \eqtag{4} we see that the choice
$$
\rho:={(\al+\be)^{\scriptstyle{3\over 2}}\over\al^\half\,\be^\half}
\eqno\eq{7}
$$
makes $C_n$ continuous on $(\al,\be)$ in the extended parameter plane.
Next we see from \eqtag{5} that the choice
$$
\si:={\al-\be\over \al+\be}
\eqno\eq{8}
$$
makes $B_n$ continuous in $(\al,\be)$ (extended) as well. Indeed,
we can now rewrite
$$
B_n={\be^{-1}-\al^{-1}\over(\al^{-1}+\be^{-1})^\half}\,
{4n+2+4n\,(n+1)/(\al+\be)\over
(1+2n/(\al+\be))\,(1+(2n+2)/(\al+\be))}\,,
$$
which is continuous in $(\al^{-1},\be^{-1})$ for $\al^{-1},\be^{-1}\ge0$.

As a result we can consider
the $(\al^{-1},\be^{-1})$-parameter plane.
For $\al^{-1},\be^{-1}>0$ we have the rescaled Jacobi polynomials
\eqtag{6} with $\rho$ and $\si$ given by \eqtag{7} and \eqtag{8}.
These polynomials extend continously to the closure
$\{(\al^{-1},\be^{-1})\mid \al^{-1},\be^{-1}\ge0\}$.
On the boundary lines $\{(\al^{-1},0)\mid\al^{-1}>0\}$
and $\{(0,\be^{-1})\mid\be^{-1}>0\}$
these polynomials become rescaled Laguerre polynomials.
On the boundary vertex $(0,0)$ the polynomials become Hermite
polynomials.

\Sec3 {Uniform limits in the Askey tableau}
The Askey tableau is given by the following chart
(cf.\ Askey \& Wilson [1, Appendix]).
$$
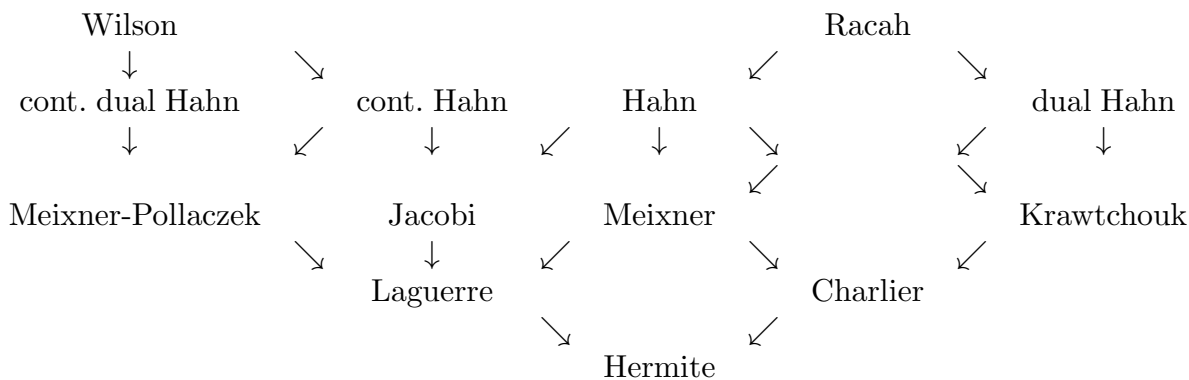
\matrix{
{\rm Wilson}&       &         &     &           &  &{\rm Racah}&& \cr
\downarrow&\searrow&	      &      &          &\swarrow& &\searrow&\cr
{\rm cont.\;dual\;Hahn}& &{\rm cont.\;Hahn}&&{\rm Hahn}&&& &{\rm dual\; Hahn}\cr
\downarrow&\swarrow&\downarrow&\swarrow&\downarrow&\searrow&&\swarrow&\downarrow
\cr
        &           &           &   &           &\swarrow&&\searrow& \cr
\hbox{ Meixner-Pollaczek}&&{\rm Jacobi}&&{\rm Meixner}& &&  &{\rm Krawtchouk}\cr
        &\searrow & \downarrow  &\swarrow&   &\searrow& &\swarrow& \cr
        &         &{\rm Laguerre}& & & &{\rm Charlier}&& \cr
   &       &        &\searrow&      &\swarrow&  &&\cr
   &       &        &       & {\rm Hermite} &    &   &   & }
$$
The various families of orthogonal polynomials mentioned here are
all of classical type, i.e., the orthogonal polynomials
$\{p_n\}_{n=0,1,\ldots}$
satisfy an equation of the form
$$
L\,p_n=\la_n\,p_n,
$$
where $L$ is some second operator (differential or difference) which does not
depend on $n$.
The arrows in the chart mean limit transitions between the various families.
The number of additional parameters on which the polynomials depend,
decreases as we go further down in the chart.
In the top row there are 4 parameters. In each subsequent row there is
one parameter less. The Hermite polynomials in the bottom row no longer
depend on parameters. The families in the left part of the chart
consist of polynomials being orthogonal \wrt an absolutely continuous
measure, while the ones in the right part are orthogonal \wrt a discrete
measure. In the case of Racah, Hahn, dual Hahn and Krawtchouk polynomials
the support
of the measure has finite cardinality, say $N+1$, and we consider
only polynomials up to degree $N$.

All the polynomials in this chart have explicit expressions as hypergeometric
functions.
For instance, {\sl Jacobi polynomials} are given by
$$
P_n^{(\al,\be)}(x)=\const\,{}_2F_1\left[{-n,n+\al+\be+1\atop\al+1};x\right].
$$
{\sl Hahn polynomials} are given by
$$
Q_n(x;\al,\be,N):=
{}_3F_2\left[{-n,n+\al+\be+1,-x\atop\al+1,-N};1\right],\quad
n=0,1,\ldots,N,
$$
while they satisfy orthogonality relations
$$
\sum_{x=0}^NQ_n(x)\,Q_m(x)\,{(\al+1)_x\,(\be+1)_{N-x}\over x!\,(N-x)!}=0,
\quad n\ne m.
$$
{\sl Racah polynomials} are given by
$$
R_n(x(x+\ga+\de+1);\al,\be,\ga,\de)=
{}_4F_3\left[{-n,n+\al+\be+1,-x,x+\ga+\de+1\atop
\al+1,\be+\de+1,\ga+1};1\right],
$$
where $\ga+1=-N$ and $n=0,1,\ldots,N$.

Now consider monic Racah polynomials
$$
r_n(x;\al,\be,-N-1,\de):=\const\,R_n(x;\al,\be,-N-1,\de)=x^n+\cdots
$$
and rescaled monic Racah polynomials
$$
p_n(x):=\rho^n\,r_n(\rho^{-1}x-\si;\al,\be,-N-1,\de).
\eqno\eq{9}
$$
Then the $p_n(x)$ satisfy recurrence relations \eqtag{3} with
$$
\eqalignno{
C_n=\rho^2\,n\,(n+\al)&\,(n+\be)\,(n+\al+\be)\,(n+\be+\de)
\cr
&\times{(\de-\al-n)\,
(n+N+\al+\be+1)\,(N+1-n)\over
(2n+\al+\be-1)\,(2n+\al+\be)^2\,(2n+\al+\be+1)}\,,&\eq{10}
\cr}
$$
$$
\eqalignno{
B_n=\rho\Bigl({(n+\al+\be+1)\,(n+\al+1)\,(n+\be+\de+1)\,(N-n)\over
(2n+\al+\be+1)\,(2n+\al+\be+2)}\qquad&
\cr
+{n\,(n+\be)\,(\de-\al-n)\,(n+N+\al+\be+1)\over(2n+\al+\be)\,(2n+\al+\be+1)}
&+\si\Bigr).&\eq{11}
\cr}
$$
Now we would like to express $\rho$ and $\si$ in such a way in terms of
$\al$, $\be$, $\de$ and $N$ that the polynomials $p_n(x)$ in \eqtag{9}
continuously
depend on these parameters up to boundaries at infinity in the
four-parameter plane, and such that Hahn polynomials and all
families which can be reached from the Hahn polynomials in the Askey tableau,
are obtained as polynomials $p_n(x)$ with parameters on the boundary.
This is a task analogous to what we did in section 2, but much more
complicated. Again, formula \eqtag{10} for $C_n$ should suggest a choice for
$\rho$ and next formula \eqtag{11} for $B_n$ should lead to the choice of $\si$.
Thus we arrive at
$$
\rho:=\Bigl({(\al+\be)^3\over\al\,\be\,(\be+\de)\,(N+\al+\be)\,(\de-\al)\,N}
\Bigr)^\half,
\eqno\eq{12}
$$
$$
\si:=-\,{N\,(\al+1)\,(\be+\de+1)\over\al+\be+2}\,.
\eqno\eq{13}
$$
Now turn from parameters $\al,\be,\de,N$ to parameters
$\al,b,d,\nu$ by the substitution
$$
\be=b\al,\quad
\de=(bd\nu+1)\,\al,\quad
N=b\,\nu.
\eqno\eq{14}
$$
Then we can prove the following theorem. Computations for this
are somewhat tedious.
Some of them I performed with the help of Maple V.

\bLP
{\bf Theorem}\quad
The rescaled monic Racah polynomials $p_n(x)$ given by \eqtag{9}, with
\eqtag{12}, \eqtag{13} and \eqtag{14} being substituted, are continuous
in $(\al^{-1},b^{-1},d^{-1},\nu^{-1})$ for
$\al^{-1},b^{-1},d^{-1},\allowbreak\nu^{-1}\ge0$.
When we restrict to any of the lower dimensional boundaries than we
obtain rescaled versions of other polynomials in the Askey tableau,
as given below.
\mPP
\halign{
\hfil#\hfil\quad&#\hfil\quad&\quad#\cr
dimension&specialization&orthogonal polynomial family\cr
\noalign{\smallskip\hrule\smallskip}
4&&Racah\cr
\noalign{\smallskip}
3&$d=\iy$&Hahn\cr
3&$\nu=\iy$&Jacobi\cr
3&$b=\iy$&Meixner\cr
3&$\al=\iy$&Krawtchouk\cr
\noalign{\smallskip}
2&$d,\nu=\iy$&Jacobi\cr
2&$d,b=\iy$&Meixner\cr
2&$d,\al=\iy$&Krawtchouk\cr
2&$\nu,b=\iy$&Laguerre\cr
2&$\nu,\al=\iy$&Hermite\cr
2&$b,\al=\iy$&Charlier\cr
\noalign{\smallskip}
1&$d,\nu,b=\iy$&Laguerre\cr
1&$d,b,\al=\iy$&Charlier\cr
1&$d,\nu,\al=\iy$&Hermite\cr
1&$\nu,b,\al=\iy$&Hermite\cr
\noalign{\smallskip}
0&$\nu,b,d,\al=\iy$&Hermite\cr}

\bPP
It is probably possible to formulate an anlogous theorem with Hahn polynomials
being replaced by dual Hahn polynomials.
We should start then with a different part of
four-parameter space for the Racah polynomials.
For Wilson polynomials we can start with three different regions
in four-parameter space. For each of these three cases there are
different limits in the Askey tableau (cf.\ [2, Table 4]).
We can hope that for each of these three cases a result analogous to
the above theorem will hold.
A further possible extension might involve $q$ as a fifth parameter.
One might also try to include the limits to non-polynomial special
functions like Bessel functions and Jacobi functions (cf.\ [2]).

\bLP
{\bf References}
\frenchspacing
\sPP
\item{[1]}
R. Askey \& J. Wilson,
{\sl Some basic hypergeometric orthogonal polynomials that generalize Jacobi
polynomials},
Memoirs Amer. Math. Soc. 54 (1985), no. 319.
\sPP
\item{[2]}
T. H. Koornwinder,
{\sl Group theoretic interpretations of Askey's scheme of hypergeometric
orthogonal polynomials},
in {\sl Orthogonal polynomials and their applications},
M. Alfaro,
J. S. Dehesa,
F. J. Marcellan,
J. L. Rubio de Francia \&
J. Vinuesa (eds.),
Lecture Notes in Math.
1329,
Springer,
1988,
pp. 46--72.
\nonfrenchspacing

\vskip 2 truecm
{\addressfont
\obeylines\parindent 5truecm
University of Amsterdam, Faculty of Mathematics and Computer Science
Plantage Muidergracht 24, 1018 TV Amsterdam, The Netherlands
email {\ttaddressfont thk@fwi.uva.nl}}

\bye